\documentclass[%
]{mpi2015-cscpreprint}

\usepackage[american]{babel}
\usepackage{graphicx}
\usepackage{amsthm}
\usepackage{damn}

\let\sidecaption\centering
\theoremstyle{definition}
\newtheorem{example}{Example}
\newtheorem{remark}[example]{Remark}

\begin{document}

\title{Diagonally-Addressed Matrix Nicknack: How to improve SpMV performance}

\author[$\dagger$]{Jens Saak \orcidlink{0000-0001-5567-9637}}
\affil[$\dagger$]{Max Planck Institute for Dynamics of Complex Technical Systems, \authorcr
  Sandtorstr. 1, 39106 Magdeburg, Germany.}

\author[$\dagger\ast$]{Jonas Schulze \orcidlink{0000-0002-2086-7686}}
\affil[$\ast$]{Corresponding author. \email{jschulze@mpi-magdeburg.mpg.de}}

\shorttitle{DAMN SpMV}
\shortauthor{J. Saak, J. Schulze}
\shortdate{}


\abstract{We suggest a technique to reduce the storage size of sparse matrices at no loss of information.
We call this technique \ac{DA} storage.
It exploits the typically low matrix bandwidth of matrices arising in applications.
For memory-bound algorithms,
this traffic reduction has direct benefits for both uni-precision and multi-precision algorithms.

In particular, we demonstrate how to apply \ac{DA} storage to the \acf{CSR} format
and compare the performance in computing the \ac{SpMV} product,
which is a basic building block of many iterative algorithms.
We investigate 1367~matrices from the SuiteSparse Matrix Collection~\cite{SuiteSparse}
fitting into the \ac{CSR} format using signed \unit{32}{\bit} indices.
More than \unit{95}{\%} of these matrices fit into the \ac{DA}-\ac{CSR} format using \unit{16}{\bit} column indices,
potentially after \ac{RCM} reordering~\cite{Cuthill1969}.
Using IEEE~754 \lstinline!double! precision scalars,
we observe a performance uplift of \unit{11}{\%} (single-threaded) or \unit{17.5}{\%} (multithreaded) on average
when the traffic exceeds the size of the last-level CPU cache.
The predicted uplift in this scenario is \unit{20}{\%}.
For traffic within the CPU's combined level 2 and level 3 caches,
the multithreaded performance uplift is over \unit{40}{\%} for a few test matrices.
}

\novelty{%
  We introduce \acf{DA} storage as a general technique to reduce sparse matrix traffic.
	Applying the \ac{DA} variant of the \ac{CSR} format to the \acf{SpMV} product,
  we demonstrate that our practical implementation achieves
  \unit{87.5}{\%} of the theoretically predicted performance uplift.}

\maketitle

\section{Introduction}

\begin{figure}
	\hfill
	\def\rainbowpoints{0/0/0, 1/2/0, 2/1/-1, 2/3/1, 3/3/0}
\begin{tikzpicture}[
        yscale=-1,
        thick,
]
        \def\lw{2pt}
        \def\radius{5pt}
        % matrix coordinates
        \begin{scope}[on background layer={
                color=gray,
                fill=gray!50,
                thick,
        }]
                \foreach \r in {0,...,3} {
                        \foreach \c in {0,...,3} {
                                \fill (\c,\r) circle [radius=0.5*\radius];
                        }
                }
        \end{scope}
        % matrix entries
        \foreach \r/\c/\yshift [count=\idx] in \rainbowpoints {
                \fill[rainbow\idx] (\c,\r) circle [radius=\radius];
        };
\end{tikzpicture}
	\hfill
	\def\rainbowpoints{0/0/0, 1/2/0, 2/1/-1, 2/3/1, 3/3/0}
\begin{tikzpicture}[
        yscale=-1,
        thick,
]
        \def\lw{2pt}
        \def\radius{5pt}
        % matrix coordinates and reference line
        \begin{scope}[on background layer={
                color=gray,
                fill=gray!50,
                thick,
        }]
                \foreach \r in {0,...,3} {
                        \foreach \c in {1,...,3} {
                                \fill (\c,\r) circle [radius=0.5*\radius];
                        }
                }
                \draw[dashed] (0,0) -- (0,3);
        \end{scope}
        % matrix entries
        \foreach \r/\c/\yshift [count=\idx] in \rainbowpoints {
                \scoped[on background layer] \draw[
                        line width=\lw,
                        rainbow\idx,
                        {Tee Barb[width=\radius,length=0pt]}-,
                        yshift = \yshift*\lw,
                ] (0,\r) -- (\c,\r);
                \fill[rainbow\idx] (\c,\r) circle [radius=\radius];
        };
\end{tikzpicture}
	\hfill
	\def\rainbowpoints{0/0/0, 1/2/0, 2/1/-1, 2/3/1, 3/3/0}
\begin{tikzpicture}[
        yscale=-1,
        thick,
]
        \def\lw{2pt}
        \def\radius{5pt}
        % matrix coordinates and reference line
        \begin{scope}[on background layer={
                color=gray,
                fill=gray!50,
                thick,
        }]
                \foreach \r in {1,...,3} {
                        \foreach \c [parse=true] in {0,...,\r-1} {
                                \fill (\c,\r) circle [radius=0.5*\radius];
                                \fill (\r,\c) circle [radius=0.5*\radius];
                        }
                }
                \draw[dashed] (0,0) -- (3,3);
        \end{scope}
        % matrix entries
        \foreach \r/\c/\yshift [count=\idx] in \rainbowpoints {
                \scoped[on background layer] \draw[
                        line width=\lw,
                        rainbow\idx,
                        {Tee Barb[width=\radius,length=0pt]}-,
                        yshift = \yshift*\lw,
                ] (\r,\r) -- (\c,\r);
                \fill[rainbow\idx] (\c,\r) circle [radius=\radius];
        };
\end{tikzpicture}
	\hspace*{\fill}
	\caption{%
		Sample matrix~(left) in \acf{CSR} storage~(middle) and \acf{DA-CSR} storage~(right).
		Colorful dots represent non-zero entries, gray dots are zero.
		Whiskers represent (column) indices with respect to a reference line~(dashed).
	}%
	\label{fig:layout}
\end{figure}
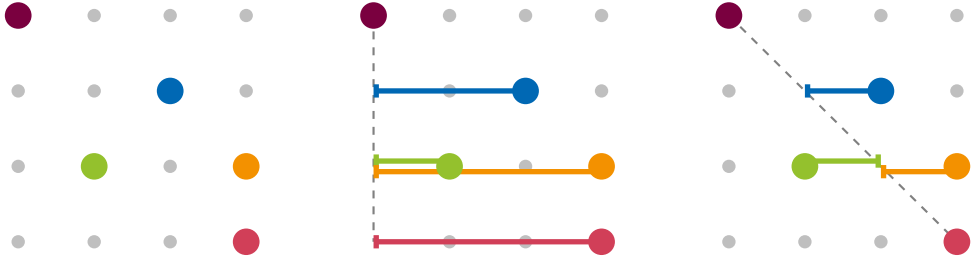
\acresetall

An operation is called memory-bound,
if its performance is limited by the memory bandwidth [Byte\per\second] of the executing hardware.
In that context,
it holds that
\begin{equation}
	\frac{\perf_1}{\perf_2} =
	\frac{\traffic_2}{\traffic_1}
	\label{eq:memory-bound}
	,
\end{equation}
where $P_i$ denotes performance [\FLOP\per\second],
and $\traffic_i$ [Byte] accounts for all the memory involved.
Hence, reducing the traffic should directly lead to a performance improvement.
It is often assumed that computing the \ac{SpMV} product
is memory-bound; see, \eg,~\cite{Goumas2008,Koza2014,Liu2013}.
The biggest contributor to the overall \ac{SpMV} traffic is the matrix.
Therefore, in the following,
we present a technique that allows for the reduction of the storage size of a sparse matrix at no loss of information.

The key ingredient of the new storage technique is the observation that
many matrices arising in, \eg, finite element simulations have a very low matrix bandwidth (under certain permutations).
That means,
potentially after permutation,
all non-zero matrix entries are located close to the matrix diagonal.
This motivates storing the indices of these entries
relative to the matrix diagonal rather than as an absolute position,
which we call \emph{\ac{DA}} storage.
Due to the small matrix bandwidth,
the relative indices may be stored in a smaller (integer) data type.
This technique is easily applicable to many sparse formats,
\eg, \acf{CSR}, \ac{BSR}, or (one of the vectors of) \ac{COO} storage.
In this paper we apply \ac{DA} storage to the \ac{CSR} format;
obtaining the \ac{DA-CSR} format;
and compare the \ac{SpMV} performance against our implementation of \ac{CSR} as well as Intel \ac{MKL}~\cite{MKL}.

\ac{DA} storage differs from Diagonal (DIA) storage
in that the new technique still requires one index per non-zero, depending on the underlying technique,
instead of one index per diagonal.
Also, it does not impose a diagonal-major order of the entries,
or require a potentially padded and full/dense storage for each of the diagonals.
The \ac{RSB} format~\cite{librsb} is another data structure for sparse matrices
motivated by a cache-efficient and parallel implementation of the \ac{SpMV} product.
It divides a sparse matrix into a tree structure of sparse blocks,
whose leaves are iterated in a Z- or Morton-order.
The leaf blocks are stored in the \ac{COO}, \ac{CSR}, or \ac{CSC} format.

\ac{RSB} was designed for arbitrary sparse matrices,
in particular, matrices without an inherent low bandwidth (under certain permutations).
\ac{RSB} allows \unit{16}{\bit} indices as well, but only for its leaf matrices.
Meanwhile, \ac{DA-CSR} has a conceptually simpler non-recursive design,
allowing \unit{16}{\bit} indices throughout,
which leads to a much lower overhead in terms of Byte per non-zero.
Therefore, \ac{DA} storage does not directly compete with the \ac{RSB} format,
but could be used in the leaf blocks within the \ac{RSB} format,
to allow for an even smaller index type.

The remainder of this paper is structured as follows.
\autoref{sec:da-storage} applies \ac{DA} storage to the \ac{CSR} format.
\autoref{sec:matrices} describes the selection of test matrices.
\autoref{sec:spmv} describes how to compute the \ac{SpMV} product using that new \ac{DA-CSR} format,
and measures the performance of \ac{SpMV}.
We conclude the paper in \autoref{sec:conclusion}.

\section{Diagonally-Addressed Storage}%
\label{sec:da-storage}

The \ac{CSR} storage of a matrix $A\in\F^{\nrows\times\ncols}$ comprises three vectors,
\cf~\autoref{snippet:csr},
\begin{lstlisting}[%
	float,
	caption={%
		(\acs{DA}-) \acs{CSR} storage of a matrix.
	},
	label=snippet:csr,
]
struct {
  oindex_t rowptr[$\nrows+1$];
  iindex_t colids[$\nnz$];
  scalar_t values[$\nnz$];
};
\end{lstlisting}
where \nnz{} denotes the number of non-zero entries,
\oindext{} and \iindext{} are integer data types,
and \scalart{} is an approximation of $\F$,
\eg, IEEE 754 \lstinline!double! or \lstinline!float! for $\F=\R$.
The \enquote{row pointers} stored in \lstinline!rowptr! and \enquote{column indices} stored in \lstinline!colids!
do the bookkeeping imposed by only storing non-zero matrix entries.

The $r$\/th entry %chktex 4
$0 \leq \text{\lstinline!rowptr[$r$]!} < \nnz $
is the index into \lstinline!colids! and \lstinline!values!
corresponding to the first non-zero of row $r$,
$0\leq r\leq\nrows$.\footnote
{The final entry \lstinline!rowptr[$\nrows$]! is set to $\nnz$ for ease of use.}
The $i$\/th entry \lstinline!colids[$i$]! is the column index
and \lstinline!values[$i$]! is the value
of the $i$\/th non-zero, $0\leq i<\nnz$.
Let $w$ denote the matrix bandwidth of $A=(a_{rc})$,
\ie~the farthest distance of a non-zero matrix entry from the matrix diagonal,
\begin{equation}
	w := \max\{ \lvert c - r\rvert : a_{rc} \neq 0\}
	\ll\ncols
	.
\end{equation}
For \ac{CSR} storage,
$\text{\lstinline!colids[$i$]!} = c_i$,
which lies in the range $[0,\ncols)$.
For \ac{DA-CSR} storage,
$\text{\lstinline!colids[$i$]!} = c_i - r_i$,
which instead lies in the range $[-w,w]$.
We illustrate this transformation in the following example.

\begin{example}
\def\rainbowcolid{\draw[line width=2pt, {Tee Barb[width=5pt,length=0pt]}-,]}
	The sample matrix shown in \autoref{fig:layout} has
	$\nrows=\ncols=4$, $\nnz=5$, and $w=1$.
	Its \ac{CSR} representation is given by
	\begin{equation}
	\left\{\begin{aligned}
		\texttt{rowptr} &= (0, 1, 2, 4, 5) \\
		\texttt{colids} &= (0, 2, 1, 3, 3)
		= (
			\tikz[baseline=0pt] \draw[line width=2pt, rainbow1] (0,0) -- (0,5pt);%
			\foreach \c [count=\i from 2] in {2, 1, 3, 3} {%
				, \tikz[scale=0.25] \rainbowcolid[rainbow\i] (0,0) -- (\c,0);
			}%
		) \\
		\texttt{values} &= (
			\tikz \fill[rainbow1] circle [radius=3pt];%
			\foreach \i in {2,...,5} {, \tikz \fill[rainbow\i] circle [radius=3pt];}%
		)
	\end{aligned}\right.
	\end{equation}
	while the \ac{DA} storage replaces \texttt{colids} to become
	\begin{equation}
		\texttt{colids}
		= (0, 1, -1, 1, 0)
		= (
			\tikz[baseline=0pt] \draw[line width=2pt, rainbow1] (0,0) -- (0,5pt);%
			\foreach \c [count=\i from 2] in {1, -1, 1} {%
				, \tikz[scale=0.25] \rainbowcolid[rainbow\i] (0,0) -- (\c,0);
			}%
			,
			\tikz[baseline=0pt] \draw[line width=2pt, rainbow5] (0,0) -- (0,5pt);%
		)
		.
	\end{equation}
	Observe that \texttt{colids} covers its full range in either storage scheme:
	$[0,\ncols-1]$ for \ac{CSR} and $[-w,w]$ for \ac{DA-CSR}.
\end{example}

\begin{figure}
	\centering
	\includegraphics[width=0.35\linewidth]{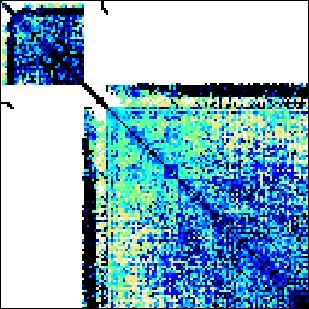}
	\hspace{0.1\linewidth}
	\includegraphics[width=0.35\linewidth]{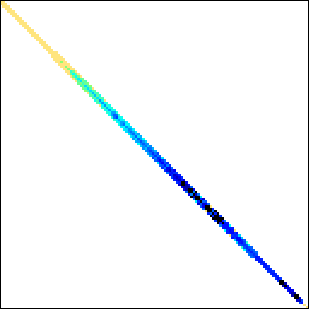}
	\caption{%
		Sparsity patterns of the matrices \texttt{GHS\_psdef/ldoor} (left)
		as well as \texttt{Janna/Bump\_2911} (right)
		from the SuiteSparse Matrix Collection~\cite{SuiteSparse}.%
	}%
	\label{fig:spy:ldoor}%
	\label{fig:spy:bump}
\end{figure}

Sometimes it is necessary to reduce the bandwidth of a matrix before \ac{DA} storage can be applied effectively,
which we observe in the next example.

%\begin{remark}[Bandwidth Reduction] % or like this?
\begin{example}
	The matrix \texttt{GHS\_posdef/ldoor} from the SuiteSparse Matrix Collection~\cite{SuiteSparse}
	has dimension 952\,203 and 46\,522\,475 pattern entries,\footnote{%
      Technically, \nnz{} refers to the number of pattern entries, which contains non-zeros as well as explicitly stored zeros.}
	distributed over a bandwidth of 686\,979.
	On average, this matrix has 49 pattern entries per row.
	See \autoref{fig:spy:ldoor} (left) for its sparsity pattern.
	Due to its block structure, the original matrix bandwidth is fairly large.
	Still, a \ac{RCM} reordering~\cite{Cuthill1969}
	reduces the bandwidth to about 9100,\footnote{%
		Our implementation yields a bandwidth of 9120,
		while the SuiteSparse Matrix Collection~\cite{SuiteSparse} reports 9134.
	}
	which is only about \unit{1}{\%} of the matrix dimension.
	Ignoring the colors, the corresponding sparsity pattern would look almost identical to \autoref{fig:spy:bump} (right).
	The reduced bandwidth is less than $2^{15} = 32\,768$ and therefore allows for the usage of \unit{16}{\bit} column indices in \ac{DA} storage,
	while both the matrix dimension (for plain \ac{CSR} storage) and the original bandwidth (for naive \ac{DA-CSR} storage) would require \unit{32}{\bit} indices.
\end{example}

Following \autoref{snippet:csr},
the matrix-related traffic amounts to
\begin{equation}
	(\nrows+1) \cdot \text{\oindexs} +
	\nnz \cdot (\text{\iindexs} + \text{\scalars})
	.
	\label{eq:traffic}
\end{equation}
If $w$ may be stored in a smaller (integer) data type than \ncols,
this allows for a smaller \iindext{} to be used.
Using an index type half the size nearly halves the bookkeeping traffic.

\begin{example}
	The matrix \texttt{Janna/Bump\_2911} from the SuiteSparse Matrix Collection~\cite{SuiteSparse} has
	dimension 2\,911\,419 and 127\,729\,899 non-zeros,
	distributed over a bandwidth of only $31\,343 < 2^{15} = 32\,768$,
	which is only about \unit{1}{\%} of the matrix dimension.
	On average, this matrix has 44 non-zeros per row.
	See \autoref{fig:spy:bump} for its sparsity pattern.
	Therefore, standard \ac{CSR} storage requires \unit{32}{\bit} indices for both
	\oindext{} and \iindext,
	which require \unit{11}{\mebi\byte} and \unit{487}{\mebi\byte} in total, respectively.
	\ac{DA-CSR} allows for \unit{16}{\bit} \iindext{} to be used,
	which requires only \unit{244}{\mebi\byte},
	thus reducing the bookkeeping traffic by $1 - \frac{11+244}{11+487} \approx \unit{48.8}{\%}$,
	or from 4.09 to \unit{2.09}{Byte} per \nnz,
	irrespective of \scalart.
	For \unit{64}{\bit} and \unit{32}{\bit} \scalart,
	\eg, IEEE~754 \lstinline!double! and \lstinline!float!,
	which in total require \unit{975}{\mebi\byte} and \unit{487}{\mebi\byte},
	using \ac{DA-CSR} instead of \ac{CSR} results in an overall matrix-related traffic reduction of
	$1 - \frac{11+244+975}{11+487+975} \approx \unit{16.5}{\%}$ and
	$1 - \frac{11+244+487}{11+487+487} \approx \unit{24.7}{\%}$, respectively.
\end{example}

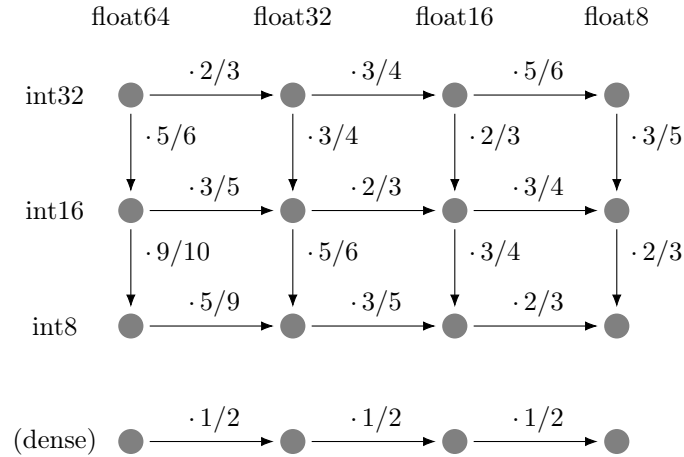
\begin{figure}
	\sidecaption
	\begin{tikzpicture}[%
	shorten >= 0.5ex,
	shorten <= 0.5ex,
	>=Latex,
]
	\matrix (grid) [%
		matrix of nodes,
		column sep = 18mm,
		row sep = 12mm,
		nodes in empty cells,
		nodes = {fill=gray,circle},
	] {
		&&& \\
		&&& \\
		&&& \\
		&&& \\
	};
	% labels: scalar_t
	\path
		(grid-1-1) node[above=8mm] {float64}
		(grid-1-2) node[above=8mm] {float32}
		(grid-1-3) node[above=8mm] {float16}
		(grid-1-4) node[above=8mm] {float8}
	;
	% labels: iindex_t
	\path
		(grid-1-1)+(-1cm,0) node {int32}
		(grid-2-1)+(-1cm,0) node {int16}
		(grid-3-1)+(-1cm,0) node {int8}
		(grid-4-1)+(-1cm,0) node {(dense)}
	;
	% iso-index
	\path[->]
		(grid-1-1) edge node[above] {$\cdot 2/3$} (grid-1-2)
		(grid-1-2) edge node[above] {$\cdot 3/4$} (grid-1-3)
		(grid-1-3) edge node[above] {$\cdot 5/6$} (grid-1-4)

		(grid-2-1) edge node[above] {$\cdot 3/5$} (grid-2-2)
		(grid-2-2) edge node[above] {$\cdot 2/3$} (grid-2-3)
		(grid-2-3) edge node[above] {$\cdot 3/4$} (grid-2-4)

		(grid-3-1) edge node[above] {$\cdot 5/9$} (grid-3-2)
		(grid-3-2) edge node[above] {$\cdot 3/5$} (grid-3-3)
		(grid-3-3) edge node[above] {$\cdot 2/3$} (grid-3-4)

		(grid-4-1) edge node[above] {$\cdot 1/2$} (grid-4-2)
		(grid-4-2) edge node[above] {$\cdot 1/2$} (grid-4-3)
		(grid-4-3) edge node[above] {$\cdot 1/2$} (grid-4-4)
	;
	% iso-scalar
	\path[->]
		(grid-1-1) edge node[pos=0.33,right] {$\cdot 5/6$} (grid-2-1)
		(grid-1-2) edge node[pos=0.33,right] {$\cdot 3/4$} (grid-2-2)
		(grid-1-3) edge node[pos=0.33,right] {$\cdot 2/3$} (grid-2-3)
		(grid-1-4) edge node[pos=0.33,right] {$\cdot 3/5$} (grid-2-4)

		(grid-2-1) edge node[pos=0.33,right] {$\cdot 9/10$} (grid-3-1)
		(grid-2-2) edge node[pos=0.33,right] {$\cdot 5/6$} (grid-3-2)
		(grid-2-3) edge node[pos=0.33,right] {$\cdot 3/4$} (grid-3-3)
		(grid-2-4) edge node[pos=0.33,right] {$\cdot 2/3$} (grid-3-4)
	;
\end{tikzpicture}
	\caption{%
		Approximate matrix-related traffic reduction when exchanging the data types used to
		store the matrix scalars~(horizontally) or column indices~(vertically)
		of (\acs{DA}-) \ac{CSR} storage,
		as well as dense storage (no indices required).
	}%
	\label{fig:traffic}
\end{figure}

For matrices with more than a few non-zeros per row, it is therefore reasonable to ignore the effect of \oindext,
\ie~to assume $\text{\oindexs}=0$.
The final percentages of the previous example would then be estimated by \(\tfrac{1}{6}\) and \(\tfrac{1}{4}\).
\autoref{fig:traffic} shows this approximate reduction in matrix-related traffic by means of formula~\eqref{eq:traffic}.
Note how smaller \iindext{};
\ie~lower bookkeeping traffic;
yield better approximations of the factor \(\tfrac{1}{2}\) observed for dense storage.\footnote{%
	Note that dense storage may be seen as having $\text{\oindexs}=\text{\iindexs}=0$,
	\ie~having zero \bit{} of bookkeeping (per \nnz).
}
Recall that by equation~\eqref{eq:memory-bound} a traffic reduction is tightly coupled with expected performance gains for memory-bound operations.

%\begin{remark}[Multi-Precision]
	While the goal of multi-precision algorithms
	is to exchange \scalart{} for a smaller data type
	(as in, \eg, \cite{Abdelfattah2021}),
	\ie~traversing the rows of \autoref{fig:traffic},
	\acs{DA} storage focuses on \iindext,
	\ie~traversing the columns of \autoref{fig:traffic}.
	However, as the main motivation of multi-precision algorithms is the memory bottleneck,
	\ac{DA} storage is expected to enable even larger speedups in that context.
	%As shown in \autoref{fig:traffic},
	%from a pure traffic perspective in the context of equation~\eqref{eq:memory-bound},
	\ac{CSR} using \unit{64}{\bit} scalars and \unit{32}{\bit} indices
	merely allows for a $\tfrac{3}{2}\times$ performance speedup when switching to \unit{32}{\bit} scalars.
	Meanwhile,
	if the matrix has a representation in \ac{DA-CSR} using \unit{16}{\bit} indices,
	the expected speedup is $\tfrac{5}{3}\times$.
	This speedup is much closer to the $2\times$ possible for dense storage, when using a scalar type half the size.
%\end{remark}

\section{Selection of Matrices}%
\label{sec:matrices}

The SuiteSparse Matrix Collection~\cite{SuiteSparse} contains
1367 square matrices having a \ac{CSR} representation using \unit{32}{\bit} indices
and a full structural rank.\footnote{
	Eventually, we are interested in using \ac{DA} storage when solving linear systems.
	We thus take full structural rank as a proxy for regularity,
	as the collection's metadata does not contain the numerical rank for all the matrices.
	Consequently, our selection of matrices contains irregular matrices as well.
}
Only 993 of these matrices (\unit{72.6}{\%}) have a dimension less than $2^{15}$,
\ie~fit into \ac{CSR} using \unit{16}{\bit} column indices.
However, for 1302 matrices (\unit{95.2}{\%}) there exists a permutation that reduces the
matrix bandwidth to below $2^{15}$,
such that these matrices fit into \ac{DA-CSR} using \unit{16}{\bit} column indices.
These are the matrices we select for further investigation.

Some of the investigated matrices are already stored in a bandwidth-reduced way.
We applied an \ac{RCM} reordering~\cite{Cuthill1969} to the ones that are not.
Unfortunately, our implementation of the \ac{RCM} permutation has not been able to
sufficiently reduce the bandwidth
of two matrices (\texttt{Janna/Long\_Coup\_dt0} and \texttt{Janna/Long\_Coup\_dt6}),
which reduces the number of matrices to 1300 (\unit{95.1}{\%}).

\section{Sparse Matrix Vector Product}%
\label{sec:spmv}

Let $A$ denote a matrix, $x$ and $y$ be vectors, and $\alpha$ and $\beta$ be scalars.
The \ac{SpMV} product denotes the operation $y \gets \alpha A x + \beta y$,
which requires
\begin{equation}
	W := 2\nnz + 2\nrows
\end{equation}
floating-point operations of \emph{work}~$W$ [\FLOP].
The \emph{performance}~$\perf$ [\FLOP\per\second] is then defined as
the ratio of work $W$ and runtime $t$,
where $t$ denotes the runtime measured in elapsed time.
The \emph{relative performance} \wrt~some baseline is computed via
\begin{equation}
	\perf\candidate/\perf\baseline
	= t\baseline/t\candidate
	,
\end{equation}
assuming $W\candidate=W\baseline$.
The \emph{traffic} [Byte] of computing the \ac{SpMV} accounts for~$x$ and~$y$
on top of the three components of $A$ in (\ac{DA}-) \ac{CSR} storage,
\cf~\autoref{snippet:csr} and formula~\eqref{eq:traffic}.
The \emph{throughput} [Byte/\second] is given by the ratio of traffic and $t$,
and the \emph{relative throughput} is then computed via
\begin{equation}
	\frac{\traffic\candidate}{\traffic\baseline}
	\cdot
	\frac{t\baseline}{t\candidate}
	,
\end{equation}
which is a scaled form of the relative performance.
Refer to \autoref{fig:traffic} for typical and approximate expected traffic ratios.
In the following, we aim to verify the predicted $\tfrac{6}{5}\times$ speedup
when replacing \unit{32}{\bit} column indices by \unit{16}{\bit} ones.

\subsection{Implementation Details and Methodology}

A prototypical implementation of the \ac{SpMV} product
for a matrix in \ac{DA-CSR} format
is given in \autoref{snippet:spmv}.
Instead of reversing the index translation in the innermost loop,
\ie~computing \lstinline!oindex_t col = row + colids[i]!,
we instead compute a shifted view~\xshift{} into the factor~\lstinline!x! one level up.
This replaces \nnz{} \oindext-additions by
\nrows{} pointer-additions.
Recall that in C/C++ the memory access
\lstinline!x[row + col]!
is equivalent to
\lstinline!*(x + (row + col))!.
Applying associativity to the computation of the pointer address,
we see that this access is also equivalent to
\lstinline!*((x + row) + col)!
and
\lstinline!xshift[col]!.

\begin{lstlisting}[float,caption={\ac{SpMV} for \ac{DA-CSR}},label=snippet:spmv]
// Input:  oindex_t *rowptr, nrows;
//         iindex_t *colids;
//         scalar_t *values, *x, *y, alpha, beta;
// Output: scalar_t *y;
for (oindex_t row = 0; row < nrows; ++row) {
  scalar_t accumulator = 0;
  scalar_t *xshift = x + row;
  for (oindex_t i = rowptr[row]; i < rowptr[row+1]; ++i) {
    iindex_t col = colids[i];
    scalar_t val = values[i];
    accumulator += val * xshift[col];
  }
  y[row] = alpha * accumulator + beta * y[row];
}
\end{lstlisting}

\begin{remark}[Non-Square Matrices]
	For tall matrices, \ie~$\nrows>\ncols$,
	\xshift{} points to memory outside \lstinline!x!,
	and must therefore never be dereferenced directly.
	Within \autoref{snippet:spmv}
	however,
	it will only be dereferenced at an offset \lstinline!col!
	that yields a memory address within \lstinline!x!.
\end{remark}

The code has been compiled with GCC~10.3.0 using \texttt{-O3 -NDEBUG -mavx2 -mfma}.
The benchmarks have been run on an Intel Xeon Skylake Silver 4110 running CentOS~7.9.2009,\footnote{\url{https://www.mpi-magdeburg.mpg.de/cluster/mechthild}}
with threads pinned using \texttt{taskset}\footnote%
{\url{https://www.man7.org/linux/man-pages/man1/taskset.1.html} and
\url{https://github.com/util-linux/util-linux/blob/master/schedutils/taskset.c}}.
Runtime measurements have been taken using nanobench~\cite{nanobench}
with \lstinline!minEpochTime! set to \unit{100}{\milli\second},
\lstinline!minEpochIterations! and \lstinline!warmup! both set to 10,
using the minimum over 11 epochs.

We measured the performance of
a naive implementation (with \nnz{} additions instead of \nrows)
as well as several implementations akin to \autoref{snippet:spmv},
optionally using OpenMP with 2, 4, 6, or 8 threads,
both for \ac{CSR} and \ac{DA-CSR} matrices.
Among all implementations executed on the given hardware,
the best performing ones were the naive implementation,
the one using 3 accumulators,
and the one using a single AVX2 accumulator (4 scalars wide) accessing \lstinline!values! using unaligned load instructions.
In the following, for each matrix, and each storage format tested,
we select the implementations and number of threads
yielding the highest performance.

%\begin{itemize}
	%\itemsep0pt
	%\item naive implementation (with \nnz{} additions instead of \nrows),
	%\item 1-8 accumulators akin to \autoref{snippet:spmv}, and
	%\item
		%vectorized accumulators initialized to either zero
		%or the first couple of summands utilizing
		%\begin{itemize}
			%\item
				%naive/unaligned loads from \lstinline!values!, or
			%\item
				%loop-peeling to be able to load data from aligned memory locations, or
			%\item
				%loop-peeling and 1-4 vector-accumulators.
		%\end{itemize}
%\end{itemize}

For the \ac{MKL}~\cite{MKL} implementation of the \ac{CSR} format,
we measured single-threaded as well as multithreaded performance.
Again, for each matrix we select the number of threads yielding the highest performance.

\subsection{Numerical Results}

For traffic within the size of the L1 cache,
a single thread yields the best performance.
Up to about \unit{100}{\kibi\byte},
which is well within the size of a single L2 cache,
the optimum number of threads increases gradually.
For traffic larger than that,
the maximum number of threads yields the best performance.
This behavior is irrespective of the matrix format and the implementation vendor
(ourselves or \ac{MKL}~\cite{MKL}).

\begin{figure}
	\includegraphics[width=\linewidth]{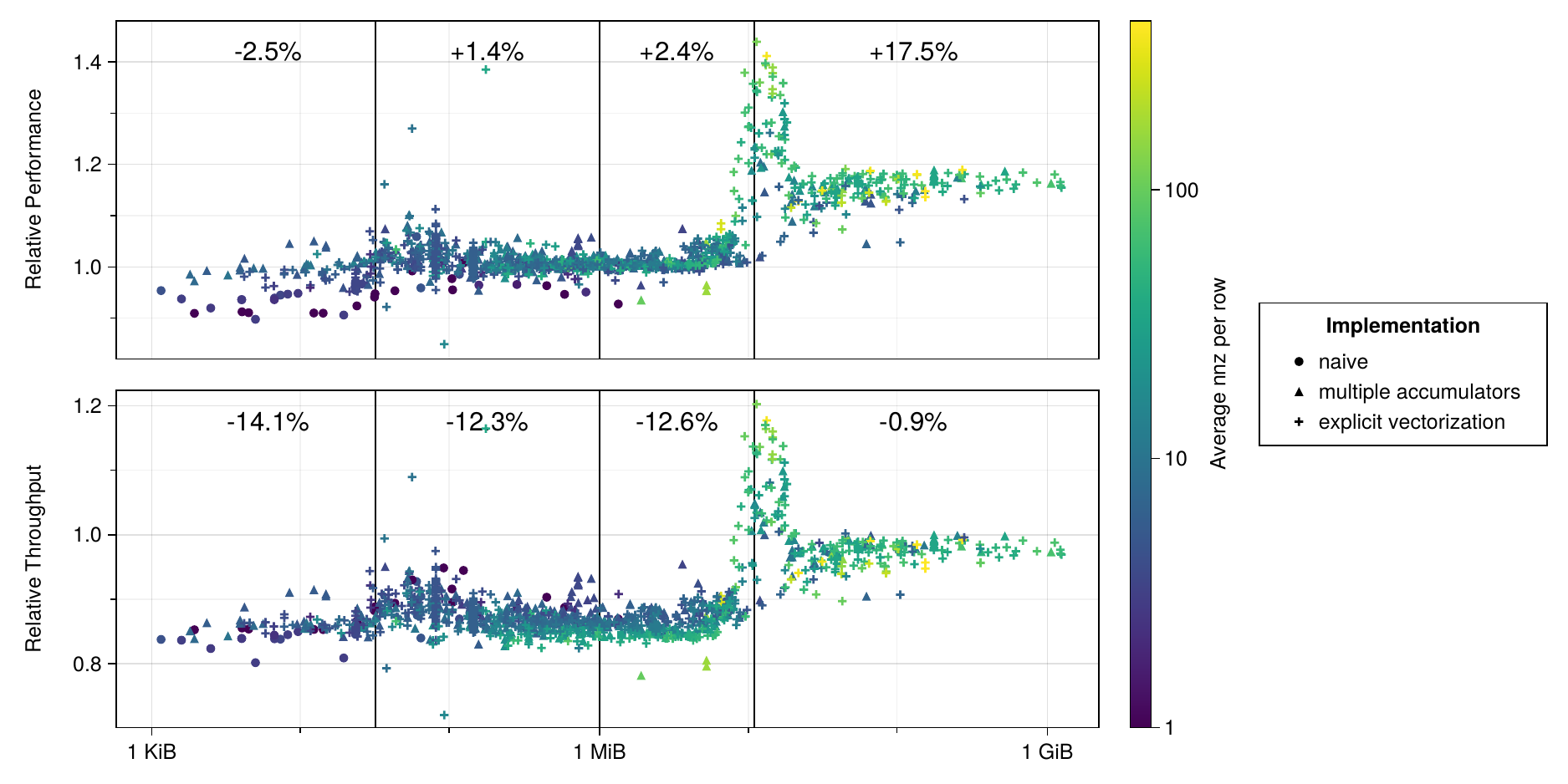}
	\caption{%
		Relative performance and throughput of \acs{SpMV} using the \acs{DA-CSR} format
		with \unit{16}{\bit} column indices
		\wrt~\ac{CSR} using \unit{32}{\bit} column indices
		as the baseline (iso-scalar).
		The sizes of the L1d, L2, and L3 CPU caches are marked with vertical lines (left to right).
	}%
	\label{fig:dacsr-vs-csr}
\end{figure}

Our implementation of the \ac{SpMV} product for the \ac{CSR} format using \unit{32}{\bit} indices
performs about the same as the \ac{MKL}~\cite{MKL},
see \autoref{tab:csr-vs-mkl}.
\autoref{fig:dacsr-vs-csr} shows the comparison of \ac{DA-CSR} using \unit{16}{\bit} column indices to \ac{CSR}.
Using \ac{DA-CSR} shows almost no change
for traffic within the combined size of the L2 caches,
\ie~up to $8 \cdot \unit{1}{\mebi\byte}$.
For traffic larger than that, up to the combined size of all caches,
\ie~up to about $\unit{8+11}{\mebi\byte}$,
we observe a larger than $1.4\times$ speedup.
For traffic beyond that,
we observe an average speedup of about \unit{+17.5}{\%},
which is reasonably close to the expected \unit{+20}{\%}.
However, the throughput drops slightly, indicating some unused potential on the given hardware.
See \autoref{tab:dacsr-vs-mkl} for the comparison of \ac{DA-CSR} to \ac{MKL}~\cite{MKL}.

\begin{table}
	\begin{minipage}{0.4\linewidth}
		\caption{%
			Average relative performance of our best \ac{SpMV} implementation
			for \ac{CSR} using \unit{32}{\bit} indices \wrt~\ac{MKL}~\cite{MKL} as the baseline.
			Values $>1$ mean we are faster.
		}
		\begin{tabular}{@{}lcc@{}}
			\hline
			Traffic & Singlethreaded & Multithreaded \\
			\hline
			%L1d   & \unit{+13.1}{\%} & \unit{+12.9}{\%} \\
			%L2    &  \unit{+7.3}{\%} &  \unit{+4.4}{\%} \\
			%L3    &  \unit{+1.2}{\%} &  \unit{+1.1}{\%} \\
			%Large &  \unit{-1.8}{\%} &  \unit{-0.6}{\%} \\
			L1d   & 1.131 & 1.129 \\
			L2    & 1.073 & 1.044 \\
			L3    & 1.012 & 1.011 \\
			Large & 0.982 & 0.994 \\
			\hline
		\end{tabular}%
		\label{tab:csr-vs-mkl}
	\end{minipage}
	\hspace{0.05\linewidth}
	\begin{minipage}{0.4\linewidth}
		\caption{%
			Average relative performance of our best \ac{SpMV} implementation
			for \ac{DA-CSR} using \unit{16}{\bit} column indices \wrt~\ac{MKL}~\cite{MKL} as the baseline.
			Values $>1$ mean we are faster.
		}
		\begin{tabular}{@{}lcc@{}}
			\hline
			Traffic & Single-threaded & Multithreaded \\
			\hline
			%L1d   & \unit{+10.3}{\%} &  \unit{+9.8}{\%} \\
			%L2    &  \unit{+7.3}{\%} &  \unit{+5.5}{\%} \\
			%L3    &  \unit{+5.2}{\%} &  \unit{+3.2}{\%} \\
			%Large & \unit{+11.0}{\%} & \unit{+17.2}{\%} \\
			L1d   & 1.103 & 1.098 \\
			L2    & 1.073 & 1.055 \\
			L3    & 1.052 & 1.032 \\
			Large & 1.110 & 1.172 \\
			\hline
		\end{tabular}%
		\label{tab:dacsr-vs-mkl}
	\end{minipage}
\end{table}

\begin{remark}[\ac{CSR} using \unit{16}{\bit} column indices]
	Recall that 933 matrices have a direct representation in \ac{CSR} using smaller column indices.
	The \ac{DA-CSR} format performs on par with \ac{CSR} using the same index types for these matrices.
\end{remark}

\section{Conclusion and Outlook}%
\label{sec:conclusion}

\acf{DA} storage allows to nearly halve the bookkeeping traffic in sparse matrix storage formats,
when the matrix bandwidth allows for an index type half the size
\unskip.
On the hardware used,
\ac{DA-CSR} storage with \unit{16}{\bit} column indices
improves the single-threaded \ac{SpMV} performance over
\ac{CSR} storage with \unit{32}{\bit} column indices
by more than \unit{17}{\%}, for both our implementation and \ac{MKL}~\cite{MKL}
if the traffic exceeds the size of the L3 cache of the CPU\@.
Meanwhile, \ac{DA-CSR} performs no worse than \ac{CSR}
when using the same data types.

\section*{Code and Data Availability}

The source code is available at:
\begin{center}
	DOI \href{https://doi.org/10.5281/zenodo.8104335}{10.5281/zenodo.8104335}
\end{center}
The visualizations in this paper have been created
using Ti\emph{k}Z~\cite{TikZ} and Makie.jl~\cite{Makie}.
The \ac{SpMV} performance measurements for the reported experiments are available at:
\begin{center}
	DOI \href{https://doi.org/10.5281/zenodo.7551699}{10.5281/zenodo.7551699}
\end{center}

\bibliographystyle{abbrv}
\IfFileExists{./DAMN-SpMV-pamm.bib}{\bibliography{DAMN-SpMV-pamm.bib}}{\bibliography{damn.bib}}

\end{document}